\documentclass{amsart}
\usepackage{amssymb}
\usepackage{amsfonts}

\newtheorem{theorem}{Theorem}
\theoremstyle{plain}

\newtheorem{corollary}{Corollary}

\newtheorem{definition}{Definition}

\newtheorem{lemma}{Lemma}

\newtheorem{remark}{Remark}

\numberwithin{equation}{section}
\input{tcilatex}

\begin{document}
\title[SOME NEW HADAMARD-TYPE INEQUALITIES]{ON\ SOME NEW HADAMARD-TYPE
INEQUALITIES FOR CO-ORDINATED QUASI-CONVEX FUNCTIONS}
\author{$^{\blacktriangledown }$M. Emin \"{O}ZDEM\.{I}R}
\address{$^{\blacktriangledown }$Atat\"{u}rk University, K.K. Education
Faculty, Department of Mathematics, 25240, Kampus, Erzurum, Turkey}
\email{emos@atauni.edu.tr}
\author{$^{\clubsuit ,\bigstar }$\c{C}etin YILDIZ}
\address{$^{\bigstar }$Atat\"{u}rk University, K.K. Education Faculty,
Department of Mathematics, 25240, Kampus, Erzurum, Turkey}
\email{yildizcetiin@yahoo.com}
\author{$^{\clubsuit }$Ahmet Ocak AKDEM\.{I}R}
\curraddr{$^{\spadesuit }$A\u{g}r\i\ \.{I}brahim \c{C}e\c{c}en University,
Faculty of Science and Letters, Department of Mathematics, 04100, A\u{g}r\i
, Turkey}
\email{ahmetakdemir@agri.edu.tr}
\subjclass[2000]{Primary 26A51, 26D15}
\keywords{Quasi-convex functions, H\"{o}lder Inequality, Power Mean
Inequality, co-ordinates, Lipschitzian function.\\
$^{\spadesuit }$Corresponding Author}

\begin{abstract}
In this paper, we establish some Hadamard-type inequalities based on
co-ordinated quasi-convexity. Also we define a new mapping associated to
co-ordinated convexity and we prove some properties of this mapping.
\end{abstract}

\maketitle

\section{INTRODUCTION}

Let $f:I\subset 
\mathbb{R}
\rightarrow 
\mathbb{R}
$ be a convex function on the interval of $I$ of real numbers and $a,b\in I$
with $a<b.$ The following double inequality%
\begin{equation}
f\left( \frac{a+b}{2}\right) \leq \frac{1}{b-a}\dint\nolimits_{a}^{b}f(x)dx%
\leq \frac{f(a)+f(b)}{2}  \label{1.1}
\end{equation}%
is well-known in the literature as Hadamard's inequality. We recall some
definitions;

\begin{definition}
(See \cite{pec}) A function $f:\left[ a,b\right] \rightarrow 
\mathbb{R}
$ is said quasi-convex on $\left[ a,b\right] $ if%
\begin{equation*}
f\left( \lambda x+(1-\lambda )y\right) \leq \max \left\{ f(x),f(y)\right\} ,%
\text{ \ \ \ \ }\left( QC\right)
\end{equation*}%
holds for all $x,y\in \left[ a,b\right] $ and $\lambda \in \lbrack 0,1].$
\end{definition}

Clearly, any convex function is quasi-convex function. Furthermore, there
exist quasi-convex functions which are not convex. In \cite{dragomir},
Dragomir defined convex functions on the co-ordinates as following:

\begin{definition}
Let us consider the bidimensional interval $\Delta =[a,b]\times \lbrack c,d]$
in $%
\mathbb{R}
^{2}$ with $a<b,$ $c<d.$ A function $f:\Delta \rightarrow 
\mathbb{R}
$ will be called convex on the co-ordinates if the partial mappings $%
f_{y}:[a,b]\rightarrow 
\mathbb{R}
,$ $f_{y}(u)=f(u,y)$ and $f_{x}:[c,d]\rightarrow 
\mathbb{R}
,$ $f_{x}(v)=f(x,v)$ are convex where defined for all $y\in \lbrack c,d]$
and $x\in \lbrack a,b].$ Recall that the mapping $f:\Delta \rightarrow 
\mathbb{R}
$ is convex on $\Delta $ if the following inequality holds, 
\begin{equation*}
f(\lambda x+(1-\lambda )z,\lambda y+(1-\lambda )w)\leq \lambda
f(x,y)+(1-\lambda )f(z,w)
\end{equation*}%
for all $(x,y),(z,w)\in \Delta $ and $\lambda \in \lbrack 0,1].$
\end{definition}

In \cite{dragomir}, Dragomir established the following inequalities of
Hadamard's type for co-ordinated convex functions on a rectangle from the
plane $%
\mathbb{R}
^{2}.$

\begin{theorem}
(see \cite{dragomir}, Theorem 1) Suppose that $f:\Delta =\left[ a,b\right]
\times \left[ c,d\right] \rightarrow 
\mathbb{R}
$ is convex on the co-ordinates on $\Delta .$ Then one has the inequalities;%
\begin{eqnarray}
f\left( \frac{a+b}{2},\frac{c+d}{2}\right) &\leq &\frac{1}{2}\left[ \frac{1}{%
b-a}\dint\nolimits_{a}^{b}f\left( x,\frac{c+d}{2}\right) dx+\frac{1}{d-c}%
\dint\nolimits_{c}^{d}f\left( \frac{a+b}{2},y\right) dy\right]  \notag \\
&\leq &\frac{1}{\left( b-a\right) \left( d-c\right) }\int_{a}^{b}%
\int_{c}^{d}f\left( x,y\right) dydx  \label{1.2} \\
&\leq &\frac{1}{4}\left[ \frac{1}{b-a}\dint\nolimits_{a}^{b}f\left(
x,c\right) dx+\frac{1}{b-a}\dint\nolimits_{a}^{b}f\left( x,d\right) dx\right.
\notag \\
&&\left. \frac{1}{d-c}\dint\nolimits_{c}^{d}f\left( a,y\right) dy+\frac{1}{%
d-c}\dint\nolimits_{c}^{d}f\left( b,y\right) dy\right]  \notag \\
&\leq &\frac{f\left( a,c\right) +f\left( b,c\right) +f\left( a,d\right)
+f\left( b,d\right) }{4}  \notag
\end{eqnarray}%
The above inequalities are sharp.
\end{theorem}

Similar results for co-ordinated $m-$convex and $(\alpha ,m)-$convex
functions can be found in \cite{OZ}. In \cite{SS}, Dragomir considered a
mapping which closely connected with above inequalities and established main
properties of this mapping as following:

Now, for a mapping $f:\Delta :=\left[ a,b\right] \times \left[ c,d\right]
\rightarrow 
\mathbb{R}
$ is convex on the co-ordinates on $\Delta ,$ we can define the mapping $H:%
\left[ 0,1\right] ^{2}\rightarrow 
\mathbb{R}
,$%
\begin{equation*}
H(t,s):=\frac{1}{\left( b-a\right) \left( d-c\right) }\dint\limits_{a}^{b}%
\dint\limits_{c}^{d}f\left( tx+(1-t)\frac{a+b}{2},sy+(1-s)\frac{c+d}{2}%
\right) dxdy
\end{equation*}

\begin{theorem}
Suppose that $f:\Delta \subset 
\mathbb{R}
^{2}\rightarrow 
\mathbb{R}
$ is convex on the co-ordinates on $\Delta =\left[ a,b\right] \times \left[
c,d\right] .$ Then:

(i) The mapping $H$ is convex on the co-ordinates on $\left[ 0,1\right]
^{2}. $

(ii) We have the bounds%
\begin{equation*}
\sup_{(t,s)\in \left[ 0,1\right] ^{2}}H(t,s)=\frac{1}{\left( b-a\right)
\left( d-c\right) }\dint\limits_{a}^{b}\dint\limits_{c}^{d}f\left(
x,y\right) dxdy=H(1,1)
\end{equation*}%
\begin{equation*}
\inf_{(t,s)\in \left[ 0,1\right] ^{2}}H(t,s)=f\left( \frac{a+b}{2},\frac{c+d%
}{2}\right) =H(0,0)
\end{equation*}

(iii) The mapping $H$ is monotonic nondecreasing on the co-ordinates.
\end{theorem}

\begin{definition}
Consider a function $f:V\rightarrow 
\mathbb{R}
$ defined on a subset $V$ of $%
\mathbb{R}
_{n}$, $n\in 
\mathbb{N}
$. Let $L=(L_{1},L_{2},...,L_{n})$ where $L_{i}\geq 0$, $i=1,2,...,n.$ We
say that $f$ is $L-$Lipschitzian function if
\end{definition}

\begin{equation*}
|f(x)-f(y)|\leq \sum\limits_{i=1}^{n}L\left\vert x_{i}-y_{i}\right\vert
\end{equation*}

for all $x,y\in V$.

In \cite{ozdemir}, \"{O}zdemir \textit{et al.} defined quasi-convex function
on the co-ordinates as following:

\begin{definition}
A function $f:\Delta =\left[ a,b\right] \times \left[ c,d\right] \rightarrow 
\mathbb{R}
$ is said quasi-convex function on the co-ordinates on $\Delta $ if the
following inequality%
\begin{equation*}
f\left( \lambda x+\left( 1-\lambda \right) z,\lambda y+\left( 1-\lambda
\right) w\right) \leq \max \left\{ f\left( x,y\right) ,f\left( z,w\right)
\right\}
\end{equation*}%
holds for all $\left( x,y\right) ,$ $\left( z,w\right) \in \Delta $ and $%
\lambda \in \left[ 0,1\right] .$
\end{definition}

Let consider a bidimensional interval $\Delta :=\left[ a,b\right] \times %
\left[ c,d\right] ,$ $f:\Delta \rightarrow 
\mathbb{R}
$ will be called co-ordinated quasi-convex on the co-ordinates if the
partial mappings%
\begin{equation*}
f_{y}:\left[ a,b\right] \rightarrow 
\mathbb{R}
,\text{ \ \ }f_{y}\left( u\right) =f\left( u,y\right)
\end{equation*}%
and%
\begin{equation*}
f_{x}:\left[ c,d\right] \rightarrow 
\mathbb{R}
,\text{ \ \ }f_{x}\left( v\right) =f\left( x,v\right)
\end{equation*}%
are convex where defined for all $y\in \left[ c,d\right] $ and $x\in \left[
a,b\right] .$ We denote by $QC(\Delta )$ the class of quasi-convex functions
on the co-ordinates on $\Delta .$

In \cite{sarikaya}, Sar\i kaya \textit{et al. }proved following Lemma and
established some inequalities for co-ordinated convex functions.

\begin{lemma}
Let $f:\Delta \subset 
\mathbb{R}
^{2}\rightarrow 
\mathbb{R}
$ be a partial differentiable mapping on $\Delta :=[a,b]\times \lbrack c,d]$
in $%
\mathbb{R}
^{2}$ with $a<b$ and $c<d.$ If $\frac{\partial ^{2}f}{\partial t\partial s}%
\in L(\Delta ),$ then the following equality holds:%
\begin{eqnarray*}
&&\frac{f\left( a,c\right) +f\left( a,d\right) +f\left( b,c\right) +f\left(
b,d\right) }{4}+\frac{1}{\left( b-a\right) \left( d-c\right) }%
\int_{a}^{b}\int_{c}^{d}f\left( x,y\right) dydx \\
&&-\frac{1}{2}\left[ \frac{1}{b-a}\dint\nolimits_{a}^{b}\left[ f\left(
x,c\right) +f\left( x,d\right) \right] dx+\frac{1}{d-c}\dint\nolimits_{c}^{d}%
\left[ f\left( a,y\right) +f\left( b,y\right) \right] dy\right] \\
&=&\frac{(b-a)(d-c)}{4}\int_{0}^{1}\int_{0}^{1}(1-2t)(1-2s)\frac{\partial
^{2}f}{\partial t\partial s}\left( ta+(1-t)b,sc+(1-s)d\right) dtds.
\end{eqnarray*}
\end{lemma}

The main purpose of this paper is to obtain some inequalities for
co-ordinated quasi-convex functions by using Lemma 1 and elemantery analysis.

\section{MAIN RESULTS}

\begin{theorem}
\label{1} Let $f:\Delta \subset 
\mathbb{R}
^{2}\rightarrow 
\mathbb{R}
$ be a partial differentiable mapping on $\Delta :=[a,b]\times \lbrack c,d]$
in $%
\mathbb{R}
^{2}$ with $a<b$ and $c<d.$ If $\left\vert \frac{\partial ^{2}f}{\partial
t\partial s}\right\vert $ is quasi-convex on the co-ordinates on $\Delta ,$
then one has the inequality:%
\begin{eqnarray*}
&&\left\vert \frac{f\left( a,c\right) +f\left( a,d\right) +f\left(
b,c\right) +f\left( b,d\right) }{4}+\frac{1}{\left( b-a\right) \left(
d-c\right) }\int_{a}^{b}\int_{c}^{d}f\left( x,y\right) dydx-A\right\vert \\
&\leq &\frac{(b-a)(d-c)}{16}\max \left\{ \left\vert \frac{\partial ^{2}f}{%
\partial t\partial s}(a,b)\right\vert ,\left\vert \frac{\partial ^{2}f}{%
\partial t\partial s}(c,d)\right\vert \right\}
\end{eqnarray*}%
where%
\begin{equation*}
A=\frac{1}{2}\left[ \frac{1}{b-a}\dint\nolimits_{a}^{b}\left[ f\left(
x,c\right) +f\left( x,d\right) \right] dx+\frac{1}{d-c}\dint\nolimits_{c}^{d}%
\left[ f\left( a,y\right) +f\left( b,y\right) \right] dy\right] .
\end{equation*}
\end{theorem}

\begin{proof}
From Lemma 1, we can write%
\begin{eqnarray*}
&&\left\vert \dfrac{f\left( a,c\right) +f\left( a,d\right) +f\left(
b,c\right) +f\left( b,d\right) }{4}\right. \\
&&\left. +\dfrac{1}{\left( b-a\right) \left( d-c\right) }\dint_{a}^{b}%
\dint_{c}^{d}f\left( x,y\right) dydx-A\right\vert \\
&\leq &\dfrac{\left( b-a\right) \left( d-c\right) }{4} \\
&&\times \dint_{0}^{1}\dint_{0}^{1}\left\vert (1-2t)(1-2s)\right\vert
\left\vert \dfrac{\partial ^{2}f}{\partial t\partial s}\left(
ta+(1-t)b,sc+(1-s)d\right) \right\vert dtds.
\end{eqnarray*}%
Since$\left\vert \frac{\partial ^{2}f}{\partial t\partial s}\right\vert $ is
quasi-convex on the co-ordinates on $\Delta $, we have%
\begin{eqnarray*}
&&\left\vert \dfrac{f\left( a,c\right) +f\left( a,d\right) +f\left(
b,c\right) +f\left( b,d\right) }{4}\right. \\
&&\left. +\dfrac{1}{\left( b-a\right) \left( d-c\right) }\dint_{a}^{b}%
\dint_{c}^{d}f\left( x,y\right) dydx-A\right\vert \\
&\leq &\dfrac{\left( b-a\right) \left( d-c\right) }{4} \\
&&\times \dint_{0}^{1}\dint_{0}^{1}\left\vert (1-2t)(1-2s)\right\vert \max
\left\{ \left\vert \frac{\partial ^{2}f}{\partial t\partial s}%
(a,b)\right\vert ,\left\vert \frac{\partial ^{2}f}{\partial t\partial s}%
(c,d)\right\vert \right\} dtds.
\end{eqnarray*}%
On the other hand, we have%
\begin{equation*}
\dint_{0}^{1}\dint_{0}^{1}\left\vert (1-2t)(1-2s)\right\vert dtds=\frac{%
(b-a)(d-c)}{16}.
\end{equation*}%
The proof is complete.
\end{proof}

\begin{theorem}
Let $f:\Delta \subset 
\mathbb{R}
^{2}\rightarrow 
\mathbb{R}
$ be a partial differentiable mapping on $\Delta :=[a,b]\times \lbrack c,d]$
in $%
\mathbb{R}
^{2}$ with $a<b$ and $c<d.$ If $\left\vert \frac{\partial ^{2}f}{\partial
t\partial s}\right\vert ^{q},$ $q>1,$ is quasi-convex function on the
co-ordinates on $\Delta ,$ then one has the inequality:%
\begin{eqnarray*}
&&\left\vert \frac{f\left( a,c\right) +f\left( a,d\right) +f\left(
b,c\right) +f\left( b,d\right) }{4}+\frac{1}{\left( b-a\right) \left(
d-c\right) }\int_{a}^{b}\int_{c}^{d}f\left( x,y\right) dydx-A\right\vert \\
&\leq &\frac{(b-a)(d-c)}{4(p+1)^{\frac{2}{p}}}\left( \max \left\{ \left\vert 
\frac{\partial ^{2}f}{\partial t\partial s}(a,b)\right\vert ^{q},\left\vert 
\frac{\partial ^{2}f}{\partial t\partial s}(c,d)\right\vert ^{q}\right\}
\right) ^{\frac{1}{q}}
\end{eqnarray*}%
where%
\begin{equation*}
A=\frac{1}{2}\left[ \frac{1}{b-a}\dint\nolimits_{a}^{b}\left[ f\left(
x,c\right) +f\left( x,d\right) \right] dx+\frac{1}{d-c}\dint\nolimits_{c}^{d}%
\left[ f\left( a,y\right) +f\left( b,y\right) \right] dy\right]
\end{equation*}%
and $\frac{1}{p}+\frac{1}{q}=1.$
\end{theorem}

\begin{proof}
From Lemma 1 and using H\"{o}lder inequality, we get%
\begin{eqnarray*}
&&\left\vert \dfrac{f\left( a,c\right) +f\left( a,d\right) +f\left(
b,c\right) +f\left( b,d\right) }{4}\right. \\
&&\left. +\dfrac{1}{\left( b-a\right) \left( d-c\right) }\dint_{a}^{b}%
\dint_{c}^{d}f\left( x,y\right) dydx-A\right\vert \\
&\leq &\dfrac{\left( b-a\right) \left( d-c\right) }{4} \\
&&\times \dint_{0}^{1}\dint_{0}^{1}\left\vert (1-2t)(1-2s)\right\vert
\left\vert \dfrac{\partial ^{2}f}{\partial t\partial s}\left(
ta+(1-t)b,sc+(1-s)d\right) \right\vert dtds \\
&\leq &\dfrac{\left( b-a\right) \left( d-c\right) }{4} \\
&&\times \left( \dint_{0}^{1}\dint_{0}^{1}\left\vert (1-2t)(1-2s)\right\vert
^{p}dtds\right) ^{\frac{1}{p}}\left( \dint_{0}^{1}\dint_{0}^{1}\left\vert 
\dfrac{\partial ^{2}f}{\partial t\partial s}\left(
ta+(1-t)b,sc+(1-s)d\right) \right\vert ^{q}dtds\right) ^{\frac{1}{q}}.
\end{eqnarray*}%
Since$\left\vert \frac{\partial ^{2}f}{\partial t\partial s}\right\vert ^{q}$
is quasi-convex on the co-ordinates on $\Delta $, we have%
\begin{eqnarray*}
&&\left\vert \dfrac{f\left( a,c\right) +f\left( a,d\right) +f\left(
b,c\right) +f\left( b,d\right) }{4}\right. \\
&&\left. +\dfrac{1}{\left( b-a\right) \left( d-c\right) }\dint_{a}^{b}%
\dint_{c}^{d}f\left( x,y\right) dydx-A\right\vert \\
&\leq &\dfrac{\left( b-a\right) \left( d-c\right) }{4} \\
&&\times \left( \dint_{0}^{1}\dint_{0}^{1}\left\vert (1-2t)(1-2s)\right\vert
^{p}dtds\right) ^{\frac{1}{p}}\left( \max \left\{ \left\vert \frac{\partial
^{2}f}{\partial t\partial s}(a,b)\right\vert ^{q},\left\vert \frac{\partial
^{2}f}{\partial t\partial s}(c,d)\right\vert ^{q}\right\} \right) ^{\frac{1}{%
q}} \\
&=&\frac{(b-a)(d-c)}{4(p+1)^{\frac{2}{p}}}\left( \max \left\{ \left\vert 
\frac{\partial ^{2}f}{\partial t\partial s}(a,b)\right\vert ^{q},\left\vert 
\frac{\partial ^{2}f}{\partial t\partial s}(c,d)\right\vert ^{q}\right\}
\right) ^{\frac{1}{q}}.
\end{eqnarray*}%
So, the proof is complete.
\end{proof}

\begin{corollary}
Since $\frac{1}{4}<\frac{1}{(p+1)^{\frac{2}{p}}}<1,$ for $p>1$, we have the
following inequality;%
\begin{eqnarray*}
&&\left\vert \frac{f\left( a,c\right) +f\left( a,d\right) +f\left(
b,c\right) +f\left( b,d\right) }{4}+\frac{1}{\left( b-a\right) \left(
d-c\right) }\int_{a}^{b}\int_{c}^{d}f\left( x,y\right) dydx-A\right\vert \\
&\leq &\frac{(b-a)(d-c)}{4}\left( \max \left\{ \left\vert \frac{\partial
^{2}f}{\partial t\partial s}(a,b)\right\vert ^{q},\left\vert \frac{\partial
^{2}f}{\partial t\partial s}(c,d)\right\vert ^{q}\right\} \right) ^{\frac{1}{%
q}}.
\end{eqnarray*}
\end{corollary}

\begin{theorem}
Let $f:\Delta \subset 
\mathbb{R}
^{2}\rightarrow 
\mathbb{R}
$ be a partial differentiable mapping on $\Delta :=[a,b]\times \lbrack c,d]$
in $%
\mathbb{R}
^{2}$ with $a<b$ and $c<d.$ If $\left\vert \frac{\partial ^{2}f}{\partial
t\partial s}\right\vert ^{q},$ $q\geq 1,$ is quasi-convex function on the
co-ordinates on $\Delta ,$ then one has the inequality:%
\begin{eqnarray*}
&&\left\vert \frac{f\left( a,c\right) +f\left( a,d\right) +f\left(
b,c\right) +f\left( b,d\right) }{4}+\frac{1}{\left( b-a\right) \left(
d-c\right) }\int_{a}^{b}\int_{c}^{d}f\left( x,y\right) dydx-A\right\vert \\
&\leq &\frac{(b-a)(d-c)}{16}\left( \max \left\{ \left\vert \frac{\partial
^{2}f}{\partial t\partial s}(a,b)\right\vert ^{q},\left\vert \frac{\partial
^{2}f}{\partial t\partial s}(c,d)\right\vert ^{q}\right\} \right) ^{\frac{1}{%
q}}
\end{eqnarray*}%
where%
\begin{equation*}
A=\frac{1}{2}\left[ \frac{1}{b-a}\dint\nolimits_{a}^{b}\left[ f\left(
x,c\right) +f\left( x,d\right) \right] dx+\frac{1}{d-c}\dint\nolimits_{c}^{d}%
\left[ f\left( a,y\right) +f\left( b,y\right) \right] dy\right] .
\end{equation*}
\end{theorem}

\begin{proof}
From Lemma 1 and using Power Mean inequality, we can write%
\begin{eqnarray*}
&&\left\vert \dfrac{f\left( a,c\right) +f\left( a,d\right) +f\left(
b,c\right) +f\left( b,d\right) }{4}\right. \\
&&\left. +\dfrac{1}{\left( b-a\right) \left( d-c\right) }\dint_{a}^{b}%
\dint_{c}^{d}f\left( x,y\right) dydx-A\right\vert \\
&\leq &\dfrac{\left( b-a\right) \left( d-c\right) }{4} \\
&&\times \dint_{0}^{1}\dint_{0}^{1}\left\vert (1-2t)(1-2s)\right\vert
\left\vert \dfrac{\partial ^{2}f}{\partial t\partial s}\left(
ta+(1-t)b,sc+(1-s)d\right) \right\vert dtds \\
&\leq &\dfrac{\left( b-a\right) \left( d-c\right) }{4}\left(
\dint_{0}^{1}\dint_{0}^{1}\left\vert (1-2t)(1-2s)\right\vert dtds\right) ^{1-%
\frac{1}{q}} \\
&&\times \left( \dint_{0}^{1}\dint_{0}^{1}\left\vert (1-2t)(1-2s)\right\vert
\left\vert \dfrac{\partial ^{2}f}{\partial t\partial s}\left(
ta+(1-t)b,sc+(1-s)d\right) \right\vert ^{q}dtds\right) ^{\frac{1}{q}}
\end{eqnarray*}%
Since$\left\vert \frac{\partial ^{2}f}{\partial t\partial s}\right\vert ^{q}$
is quasi-convex on the co-ordinates on $\Delta $, we have%
\begin{eqnarray*}
&&\left\vert \dfrac{f\left( a,c\right) +f\left( a,d\right) +f\left(
b,c\right) +f\left( b,d\right) }{4}\right. \\
&&\left. +\dfrac{1}{\left( b-a\right) \left( d-c\right) }\dint_{a}^{b}%
\dint_{c}^{d}f\left( x,y\right) dydx-A\right\vert \\
&\leq &\dfrac{\left( b-a\right) \left( d-c\right) }{4}\left( \max \left\{
\left\vert \frac{\partial ^{2}f}{\partial t\partial s}(a,b)\right\vert
^{q},\left\vert \frac{\partial ^{2}f}{\partial t\partial s}(c,d)\right\vert
^{q}\right\} \right) ^{\frac{1}{q}} \\
&&\times \left( \dint_{0}^{1}\dint_{0}^{1}\left\vert (1-2t)(1-2s)\right\vert
dtds\right) ^{1-\frac{1}{q}}\left( \dint_{0}^{1}\dint_{0}^{1}\left\vert
(1-2t)(1-2s)\right\vert dtds\right) ^{\frac{1}{q}} \\
&=&\frac{(b-a)(d-c)}{16}\left( \max \left\{ \left\vert \frac{\partial ^{2}f}{%
\partial t\partial s}(a,b)\right\vert ^{q},\left\vert \frac{\partial ^{2}f}{%
\partial t\partial s}(c,d)\right\vert ^{q}\right\} \right) ^{\frac{1}{q}}
\end{eqnarray*}%
which completes the proof.
\end{proof}

\begin{remark}
Since $\frac{1}{4}<\frac{1}{(p+1)^{\frac{2}{p}}}<1,$ for $p>1$, the
estimation in Theorem 4 is better than Theorem 3.
\end{remark}

Now, for a mapping $f:\Delta :=\left[ a,b\right] \times \left[ c,d\right]
\rightarrow 
\mathbb{R}
$ is convex on the co-ordinates on $\Delta ,$ we can define the mapping $G:%
\left[ 0,1\right] ^{2}\rightarrow 
\mathbb{R}
,$%
\begin{eqnarray*}
G(t,s) &:&=\frac{1}{4}\left[ f\left( ta+(1-t)\frac{a+b}{2},sc+(1-s)\frac{c+d%
}{2}\right) \right. \\
&&+f\left( tb+(1-t)\frac{a+b}{2},sc+(1-s)\frac{c+d}{2}\right) \\
&&+f\left( ta+(1-t)\frac{a+b}{2},sd+(1-s)\frac{c+d}{2}\right) \\
&&\left. +f\left( tb+(1-t)\frac{a+b}{2},sd+(1-s)\frac{c+d}{2}\right) \right]
\end{eqnarray*}%
We will give following theorem which contains some properties of this
mapping.

\begin{theorem}
Suppose that $f:\Delta \subset 
\mathbb{R}
^{2}\rightarrow 
\mathbb{R}
$ is convex on the co-ordinates on $\Delta =\left[ a,b\right] \times \left[
c,d\right] .$ Then:

(i) The mapping $G$ is convex on the co-ordinates on $\left[ 0,1\right]
^{2}. $

(ii) We have the bounds%
\begin{equation*}
\inf_{(t,s)\in \left[ 0,1\right] ^{2}}G(t,s)=f\left( \frac{a+b}{2},\frac{c+d%
}{2}\right) =G(0,0)
\end{equation*}%
\begin{equation*}
\sup_{(t,s)\in \left[ 0,1\right] ^{2}}G(t,s)=\dfrac{f\left( a,c\right)
+f\left( a,d\right) +f\left( b,c\right) +f\left( b,d\right) }{4}=G(1,1)
\end{equation*}

(iii) If $f$ is satisfy Lipschitzian conditions, then the mapping $G$ is $L-$%
Lipschitzian on $[0,1]\times \lbrack 0,1].$

(iv) Following inequality holds;%
\begin{eqnarray*}
&&\frac{1}{(b-a)(d-c)}\int\limits_{a}^{b}\int\limits_{c}^{d}f(x,y)dydx \\
&\leq &\frac{1}{4}\left[ \frac{f\left( a,c\right) +f\left( b,c\right)
+f(a,d)+f(b,d)}{4}\right.  \\
&&\left. +\frac{f\left( \frac{a+b}{2},c\right) +f\left( \frac{a+b}{2}%
,d\right) }{2}+f\left( \frac{a+b}{2},\frac{c+d}{2}\right) \right] .
\end{eqnarray*}
\end{theorem}

\begin{proof}
(i) Let $s\in \left[ 0,1\right] .$ For all $\alpha ,\beta \geq 0$ with $%
\alpha +\beta =1$ and $t_{1},t_{2}\in \left[ 0,1\right] ,$ then we have%
\begin{eqnarray*}
&&G(\alpha t_{1}+\beta t_{2},s) \\
&=&\frac{1}{4}\left[ f\left( \left( \alpha t_{1}+\beta t_{2}\right)
a+(1-\left( \alpha t_{1}+\beta t_{2}\right) )\frac{a+b}{2},sc+(1-s)\frac{c+d%
}{2}\right) \right. \\
&&+f\left( \left( \alpha t_{1}+\beta t_{2}\right) b+(1-\left( \alpha
t_{1}+\beta t_{2}\right) )\frac{a+b}{2},sc+(1-s)\frac{c+d}{2}\right) \\
&&+f\left( \left( \alpha t_{1}+\beta t_{2}\right) a+(1-\left( \alpha
t_{1}+\beta t_{2}\right) )\frac{a+b}{2},sd+(1-s)\frac{c+d}{2}\right) \\
&&\left. +f\left( \left( \alpha t_{1}+\beta t_{2}\right) b+(1-\left( \alpha
t_{1}+\beta t_{2}\right) )\frac{a+b}{2},sd+(1-s)\frac{c+d}{2}\right) \right]
\\
&=&\frac{1}{4}\left[ f\left( \alpha \left( t_{1}a+(1-t_{1})\frac{a+b}{2}%
\right) +\beta \left( t_{2}a+(1-t_{2})\frac{a+b}{2}\right) ,sc+(1-s)\frac{c+d%
}{2}\right) \right. \\
&&+f\left( \alpha \left( t_{1}b+(1-t_{1})\frac{a+b}{2}\right) +\beta \left(
t_{2}b+(1-t_{2})\frac{a+b}{2}\right) ,sc+(1-s)\frac{c+d}{2}\right) \\
&&+f\left( \alpha \left( t_{1}a+(1-t_{1})\frac{a+b}{2}\right) +\beta \left(
t_{2}a+(1-t_{2})\frac{a+b}{2}\right) ,sd+(1-s)\frac{c+d}{2}\right) \\
&&\left. +f\left( \alpha \left( t_{1}b+(1-t_{1})\frac{a+b}{2}\right) +\beta
\left( t_{2}b+(1-t_{2})\frac{a+b}{2}\right) ,sd+(1-s)\frac{c+d}{2}\right) %
\right] .
\end{eqnarray*}%
Using the convexity of $f,$ we obtain%
\begin{eqnarray*}
G(\alpha t_{1}+\beta t_{2},s) &\leq &\frac{1}{4}\left[ \alpha \left( f\left(
t_{1}a+(1-t_{1})\frac{a+b}{2},sc+(1-s)\frac{c+d}{2}\right) \right. \right. \\
&&+f\left( t_{1}b+(1-t_{1})\frac{a+b}{2},sc+(1-s)\frac{c+d}{2}\right) \\
&&+f\left( t_{1}a+(1-t_{1})\frac{a+b}{2},sd+(1-s)\frac{c+d}{2}\right) \\
&&\left. +f\left( t_{1}b+(1-t_{1})\frac{a+b}{2},sd+(1-s)\frac{c+d}{2}\right)
\right) \\
&&+\beta \left( f\left( t_{2}a+(1-t_{2})\frac{a+b}{2},sc+(1-s)\frac{c+d}{2}%
\right) \right. \\
&&+f\left( t_{2}b+(1-t_{2})\frac{a+b}{2},sc+(1-s)\frac{c+d}{2}\right) \\
&&+f\left( t_{2}a+(1-t_{2})\frac{a+b}{2},sd+(1-s)\frac{c+d}{2}\right) \\
&&\left. \left. +f\left( t_{2}b+(1-t_{2})\frac{a+b}{2},sd+(1-s)\frac{c+d}{2}%
\right) \right) \right] \\
&=&\alpha G(t_{1},s)+\beta G(t_{2},s).
\end{eqnarray*}%
If $s\in \left[ 0,1\right] .$ For all $\alpha ,\beta \geq 0$ with $\alpha
+\beta =1$ and $t_{1},t_{2}\in \left[ 0,1\right] ,$ then we also have;%
\begin{equation*}
G(t,\alpha s_{1}+\beta s_{2})\leq \alpha G(t,s_{1})+\beta G(t,s_{2})
\end{equation*}

and the statement is proved.

(ii) It is easy to see that by taking $t=s=0$ and $t=s=1,$ respectively, in $%
G,$ we have the bounds%
\begin{equation*}
\inf_{(t,s)\in \left[ 0,1\right] ^{2}}G(t,s)=f\left( \frac{a+b}{2},\frac{c+d%
}{2}\right) =G(0,0)
\end{equation*}%
\begin{equation*}
\sup_{(t,s)\in \left[ 0,1\right] ^{2}}G(t,s)=\dfrac{f\left( a,c\right)
+f\left( a,d\right) +f\left( b,c\right) +f\left( b,d\right) }{4}=G(1,1).
\end{equation*}%
(iii) Let $t_{1},t_{2},s_{1},s_{2}\in \left[ 0,1\right] ,$ then we have%
\begin{eqnarray*}
&&|G(t_{2},s_{2})-G(t_{1},s_{1})| \\
&=&\frac{1}{4}\left\vert f\left( t_{2}a+(1-t_{2})\frac{a+b}{2}%
,s_{2}c+(1-s_{2})\frac{c+d}{2}\right) +f\left( t_{2}b+(1-t_{2})\frac{a+b}{2}%
,s_{2}c+(1-s_{2})\frac{c+d}{2}\right) \right\vert \\
&&+f\left( t_{2}a+(1-t_{2})\frac{a+b}{2},s_{2}d+(1-s_{2})\frac{c+d}{2}%
\right) +f\left( t_{2}b+(1-t_{2})\frac{a+b}{2},s_{2}d+(1-s_{2})\frac{c+d}{2}%
\right) \\
&&-f\left( t_{1}a+(1-t_{1})\frac{a+b}{2},s_{1}c+(1-s_{1})\frac{c+d}{2}%
\right) -f\left( t_{1}b+(1-t_{1})\frac{a+b}{2},s_{1}c+(1-s_{1})\frac{c+d}{2}%
\right) \\
&&\left. -f\left( t_{1}a+(1-t_{1})\frac{a+b}{2},s_{1}d+(1-s_{1})\frac{c+d}{2}%
\right) -f\left( t_{1}b+(1-t_{1})\frac{a+b}{2},s_{1}d+(1-s_{1})\frac{c+d}{2}%
\right) \right\vert .
\end{eqnarray*}%
By using the triangle inequality, we get%
\begin{eqnarray*}
&&|G(t_{2},s_{2})-G(t_{1},s_{1})| \\
&\leq &\frac{1}{4}\left\vert f\left( t_{2}a+(1-t_{2})\frac{a+b}{2}%
,s_{2}c+(1-s_{2})\frac{c+d}{2}\right) -f\left( t_{1}a+(1-t_{1})\frac{a+b}{2}%
,s_{1}c+(1-s_{1})\frac{c+d}{2}\right) \right\vert \\
&&\left\vert +f\left( t_{2}b+(1-t_{2})\frac{a+b}{2},s_{2}c+(1-s_{2})\frac{c+d%
}{2}\right) -f\left( t_{1}b+(1-t_{1})\frac{a+b}{2},s_{1}c+(1-s_{1})\frac{c+d%
}{2}\right) \right\vert \\
&&\left\vert +f\left( t_{2}a+(1-t_{2})\frac{a+b}{2},s_{2}d+(1-s_{2})\frac{c+d%
}{2}\right) -f\left( t_{1}a+(1-t_{1})\frac{a+b}{2},s_{1}d+(1-s_{1})\frac{c+d%
}{2}\right) \right\vert \\
&&\left\vert +f\left( t_{2}b+(1-t_{2})\frac{a+b}{2},s_{2}d+(1-s_{2})\frac{c+d%
}{2}\right) -f\left( t_{1}b+(1-t_{1})\frac{a+b}{2},s_{1}d+(1-s_{1})\frac{c+d%
}{2}\right) \right\vert .
\end{eqnarray*}%
By using the $f$ is satisfy Lipschitzian conditions, then we obtain%
\begin{eqnarray*}
&&\frac{1}{4}\left\vert f\left( t_{2}a+(1-t_{2})\frac{a+b}{2}%
,s_{2}c+(1-s_{2})\frac{c+d}{2}\right) -f\left( t_{1}a+(1-t_{1})\frac{a+b}{2}%
,s_{1}c+(1-s_{1})\frac{c+d}{2}\right) \right\vert \\
&&\left\vert +f\left( t_{2}b+(1-t_{2})\frac{a+b}{2},s_{2}c+(1-s_{2})\frac{c+d%
}{2}\right) -f\left( t_{1}b+(1-t_{1})\frac{a+b}{2},s_{1}c+(1-s_{1})\frac{c+d%
}{2}\right) \right\vert \\
&&\left\vert +f\left( t_{2}a+(1-t_{2})\frac{a+b}{2},s_{2}d+(1-s_{2})\frac{c+d%
}{2}\right) -f\left( t_{1}a+(1-t_{1})\frac{a+b}{2},s_{1}d+(1-s_{1})\frac{c+d%
}{2}\right) \right\vert \\
&&\left\vert +f\left( t_{2}b+(1-t_{2})\frac{a+b}{2},s_{2}d+(1-s_{2})\frac{c+d%
}{2}\right) -f\left( t_{1}b+(1-t_{1})\frac{a+b}{2},s_{1}d+(1-s_{1})\frac{c+d%
}{2}\right) \right\vert \\
&\leq &\frac{1}{4}\left[ L_{1}(b-a)\left\vert t_{2}-t_{1}\right\vert
+L_{2}(d-c)\left\vert s_{2}-s_{1}\right\vert +L_{3}(b-a)\left\vert
t_{2}-t_{1}\right\vert +L_{4}(d-c)\left\vert s_{2}-s_{1}\right\vert \right.
\\
&&\left. +L_{5}(b-a)\left\vert t_{2}-t_{1}\right\vert +L_{6}(d-c)\left\vert
s_{2}-s_{1}\right\vert +L_{7}(b-a)\left\vert t_{2}-t_{1}\right\vert
+L_{8}(d-c)\left\vert s_{2}-s_{1}\right\vert \right] \\
&=&\frac{1}{4}\left[ \left( L_{1}+L_{2}+L_{3}+L_{4}\right) (b-a)\left\vert
t_{2}-t_{1}\right\vert +\left( L_{5}+L_{6}+L_{7}+L_{8}\right)
(d-c)\left\vert s_{2}-s_{1}\right\vert \right]
\end{eqnarray*}%
this imply that the mapping $G$ is $L-$Lipschitzian on $[0,1]\times \lbrack
0,1].$

(iv) By using the convexity of $G$ on $[0,1]\times \lbrack 0,1],$ we have%
\begin{eqnarray*}
&&f\left( ta+(1-t)\frac{a+b}{2},sc+(1-s)\frac{c+d}{2}\right) +f\left(
tb+(1-t)\frac{a+b}{2},sc+(1-s)\frac{c+d}{2}\right)  \\
&&\left. +f\left( ta+(1-t)\frac{a+b}{2},sd+(1-s)\frac{c+d}{2}\right)
+f\left( tb+(1-t)\frac{a+b}{2},sd+(1-s)\frac{c+d}{2}\right) \right]  \\
&\leq &tsf\left( a,c\right) +t(1-s)f\left( a,\frac{c+d}{2}\right)
+(1-t)sf\left( \frac{a+b}{2},c\right) +(1-t)(1-s)f\left( \frac{a+b}{2},\frac{%
c+d}{2}\right)  \\
&&+tsf\left( b,c\right) +t(1-s)f\left( b,\frac{c+d}{2}\right) +(1-t)sf\left( 
\frac{a+b}{2},c\right) +(1-t)(1-s)f\left( \frac{a+b}{2},\frac{c+d}{2}\right) 
\\
&&+tsf\left( a,d\right) +t(1-s)f\left( a,\frac{c+d}{2}\right) +(1-t)sf\left( 
\frac{a+b}{2},d\right) +(1-t)(1-s)f\left( \frac{a+b}{2},\frac{c+d}{2}\right) 
\\
&&+tsf\left( b,d\right) +t(1-s)f\left( b,\frac{c+d}{2}\right) +(1-t)sf\left( 
\frac{a+b}{2},d\right) +(1-t)(1-s)f\left( \frac{a+b}{2},\frac{c+d}{2}\right)
.
\end{eqnarray*}%
By integrating both sides of the above inequality and by taking into account
the change of the variables, we obtain 
\begin{eqnarray*}
&&\frac{1}{(b-a)(d-c)}\int\limits_{a}^{b}\int\limits_{c}^{d}f(x,y)dydx \\
&\leq &\frac{1}{4}\left[ \frac{f\left( a,c\right) +f\left( b,c\right)
+f(a,d)+f(b,d)}{4}\right.  \\
&&\left. +\frac{f\left( \frac{a+b}{2},c\right) +f\left( \frac{a+b}{2}%
,d\right) }{2}+f\left( \frac{a+b}{2},\frac{c+d}{2}\right) \right] .
\end{eqnarray*}%
Which completes the proof.
\end{proof}

\end{document}